\documentclass{gtart}

\def\ifplaintex{\expandafter\ifx\csname documentclass\endcsname\relax}

\ifplaintex 
\else
\fi

\expandafter\ifx\csname beginpicture\endcsname\relax
\expandafter\ifx\csname documentclass\endcsname\relax
\input pictex \else
\input prepictex \input pictex \input postpictex \fi\fi

\def\gt{{\mathsurround=0pt\it $\cal G\mskip-2mu$eometry \&\ 
$\cal T\!\!$opology}}        

\def\gtp{{\mathsurround=0pt\it $\cal G\mskip-2mu$eometry \&\ 
$\cal T\!\!$opology $\cal P\!$ublications}}  

\def\volumenumber#1{\def\thevolumenumber{#1}}
\def\papernumber#1{\def\thepapernumber{#1}}
\def\volumeyear#1{\def\thevolumeyear{#1}}

\def\pagenumbers#1#2{\def\startpage{#1}\def\finishpage{#2}}
\def\published#1{\def\publishdate{#1}}
\def\proposed#1{\def\theproposer{#1}}
\def\seconded#1{\def\theseconders{#1}}
\def\received#1{\def\receiveddate{#1}}

\def\accepted#1{\def\accepteddate{#1}}
\def\asciititle#1{\def\theasciititle{#1}}
\def\covertitle#1{\def\thecovertitle{#1}}

\def\asciiaddress#1{\def\theasciiaddress{#1}}

\long\def\asciiabstract#1{\long\def\theasciiabstract{#1}}
\def\asciikeywords#1{\def\theasciikeywords{#1}}

\def\shorttitle#1{\def\theshorttitle{#1}}

\let\\\par
\let\thevolumenumber\relax\let\thepapernumber\relax
\let\thevolumeyear\relax\let\thesamplenumber\relax\let\startpage\relax
\let\finishpage\relax\let\publishdate\relax\let\receiveddate\relax
\let\reviseddate\relax\let\accepteddate\relax\let\theasciititle\relax
\let\thecovertitle\relax\let\theasciiauthors\relax\let\theasciiaddress\relax
\let\theasciiabstract\relax\let\theasciikeywords\relax
\let\theasciiemail\relax\let\theshortauthors\relax\let\theshorttitle\relax

\long\def\maketitlep{   

\count0=\startpage

\gt\hfill      
\beginpicture
\setcoordinatesystem units <0.33truein, 0.33truein> point at 2.2 0.9
\setplotsymbol ({$\cal G$})
\plotsymbolspacing=9truept
\circulararc 315 degrees from 0 1 center at 0 0
\setplotsymbol ({$\cal T$})
\circulararc 315 degrees from 1 -1 center at 1 0
\endpicture
\break
{\small\ifx\thesamplenumber\relax 
Volume \else Sample
\fi\thevolumenumber\ (\thevolumeyear)
\startpage--\finishpage\nl
Published: \publishdate}
\vglue 0.5truein plus 0.4fil minus 0.1truein

{\parskip=0pt\leftskip 0pt plus 1fil\def\\{\par\smallskip}{\ifplaintex\large
\else\Large\fi\bf\thetitle}\par\medskip}   

\vglue 0pt plus 0.1fil 

{\parskip=0pt\leftskip 0pt plus 1fil\def\\{\par}{\sc\theauthors}
\par\medskip}

\vglue 0pt plus 0.1fil 

{\small\parskip=0pt\let\newline\\
{\leftskip 0pt plus 1fil\def\\{\par}{\sl\theaddress}\par}
\expandafter\ifx\theemail\relax    
\relax\else\vglue 5pt plus 0.02fil minus 2pt\def\\{\stdspace{\rm 
and}\stdspace} 
\cl{Email:\stdspace\tt\theemail}\fi
\ifx\theurl\relax                  
\relax\else\vglue 5pt plus 0.02fil minus 2pt\def\\{\stdspace{\rm 
and}\stdspace}
\cl{URL:\stdspace\tt\theurl}\fi\par}

\vglue 7pt plus 0.3fil minus 3pt

{\bf Abstract}
\vglue 5pt plus 0.1fil minus 2pt

\theabstract

\vglue 7pt plus 0.3fil minus 3pt

{\bf AMS Classification numbers}\quad Primary:\quad \theprimaryclass

Secondary:\quad \thesecondaryclass

\vglue 5pt plus 0.3fil minus 2pt

{\bf Keywords:}\quad \thekeywords

\vglue 10pt plus 0.5fil minus 5pt

{\small  Proposed: \theproposer\hfill Received: \receiveddate\nl
Seconded: \theseconders\hfill 
\ifx\reviseddate\relax                         
Accepted: \accepteddate                        
\else
Revised: \reviseddate                          
\fi}
\eject
}       

\let\maketitlepage\maketitlep
\let\maketitle\maketitlepage

\font\phead=cmsl9 scaled 950
\font=cmsl9 scaled 1050
\font\pnum=cmbx10 scaled 913
\font\lnum=cmbx10 
\font\pfoot=cmsl9 scaled 950
\font=cmsl9 scaled 1050
\ifplaintex
\headline{\vbox to 0pt{\vskip -4.5mm\line{\small\phead\ifnum
\count0=\startpage ISSN 1364-0380 (on line)
1465-3060 (printed) \hfill {\pnum\folio}\else\ifodd\count0\def\\{ }%
\ifx\theshorttitle\relax\thetitle\else\theshorttitle\fi\hfill{\pnum\folio}
\else\def\\{ and }{\pnum\folio}\hfill\ifx\theshortauthors\relax\theauthors
\else\theshortauthors\fi\fi\fi}\vss}}
\footline{\vbox to 0pt{\vglue 0mm\line{\small\pfoot\ifnum\count0=\startpage
\copyright\ \gtp\hfill\else
\gt, Volume \thevolumenumber\ (\thevolumeyear)\hfill\fi}\vss
}}
\else
\makeatletter
\def\@oddhead{{\small\lhead\ifnum\count0=\startpage ISSN 1364-0380 (on line)
1465-3060 (printed) \hfill {\lnum\number\count0}\else\ifodd\count0
\def\\{ }\ifx\theshorttitle\relax \thetitle \else\theshorttitle\fi\hfill
{\lnum\number\count0}\else\def\\{ and }{\lnum\number\count0}
\hfill\ifx\theshortauthors\relax 
\theauthors\else\theshortauthors\fi\fi\fi}}\def\@evenhead{\@oddhead}
\def\@oddfoot{\small\lfoot\ifnum\count0=\startpage\copyright\ \gtp\hfill\else
\gt, Volume \thevolumenumber\ (\thevolumeyear)\hfill\fi}
\def\@evenfoot{\@oddfoot}
\makeatother
\fi

\newwrite\gtoutfile
\long\gdef\makeheadfile{  
{\def\\{, }\def\s{ }
\immediate\openout\gtoutfile head.xxx
\immediate\write\gtoutfile{To: math@arxiv.org}
\immediate\write\gtoutfile{Subject: put or rep NNNNN:pppp}
\immediate\write\gtoutfile{--text follows this line--}
\immediate\write\gtoutfile{Proxy-for: \ifx\theasciiauthors\relax
\theauthors\else\theasciiauthors\fi\s<\ifx\theasciiemail\relax\theemail\else\theasciiemail\fi>}
\immediate\write\gtoutfile{\noexpand\\}
\immediate\write\gtoutfile{Authors: \ifx\theasciiauthors\relax
\theauthors\else\theasciiauthors\fi}
\immediate\write\gtoutfile{Title: \ifx\theasciititle\relax
\thetitle\else\theasciititle\fi}
\immediate\write\gtoutfile{Subj-class: GT or SG or MG etc}
\immediate\write\gtoutfile{MSC-class: \theprimaryclass\ifx\thesecondaryclass\relax\else, \thesecondaryclass\fi}
\immediate\write\gtoutfile{Journal-ref: Geom. Topol. \thevolumenumber
(\thevolumeyear) \startpage-\finishpage}
\immediate\write\gtoutfile{Comments: Published in Geometry and Topology at}
\immediate\write\gtoutfile{\s\s http://www.maths.warwick.ac.uk/gt/GTVol\thevolumenumber/paper\thepapernumber.abs.html}
\immediate\write\gtoutfile{\noexpand\\}
\immediate\write\gtoutfile{}
\ifx\theasciiabstract\relax
\immediate\write\gtoutfile{\theabstract}\else
\immediate\write\gtoutfile{\theasciiabstract}\fi
\immediate\write\gtoutfile{}
\immediate\write\gtoutfile{\noexpand\\}
\immediate\write\gtoutfile{}
\immediate\write\gtoutfile{<uuencoded .tar.gz file here>}
\immediate\write\gtoutfile{}
\immediate\closeout\gtoutfile}}  

\def\maketitlepage{\maketitlep\makeheadfile}
\let\maketitle\maketitlepage

\def\ifplaintex{\expandafter\ifx\csname documentclass\endcsname\relax}

\ifplaintex 
\else
\fi

\expandafter\ifx\csname beginpicture\endcsname\relax
\expandafter\ifx\csname documentclass\endcsname\relax
\input pictex \else
\input prepictex \input pictex \input postpictex \fi\fi

\def\gt{{\mathsurround=0pt\it $\cal G\mskip-2mu$eometry \&\ 
$\cal T\!\!$opology}}        

\def\gtp{{\mathsurround=0pt\it $\cal G\mskip-2mu$eometry \&\ 
$\cal T\!\!$opology $\cal P\!$ublications}}  

\def\volumenumber#1{\def\thevolumenumber{#1}}
\def\papernumber#1{\def\thepapernumber{#1}}
\def\volumeyear#1{\def\thevolumeyear{#1}}

\def\pagenumbers#1#2{\def\startpage{#1}\def\finishpage{#2}}
\def\published#1{\def\publishdate{#1}}
\def\proposed#1{\def\theproposer{#1}}
\def\seconded#1{\def\theseconders{#1}}
\def\received#1{\def\receiveddate{#1}}

\def\accepted#1{\def\accepteddate{#1}}
\def\asciititle#1{\def\theasciititle{#1}}
\def\covertitle#1{\def\thecovertitle{#1}}

\def\asciiaddress#1{\def\theasciiaddress{#1}}

\long\def\asciiabstract#1{\long\def\theasciiabstract{#1}}
\def\asciikeywords#1{\def\theasciikeywords{#1}}

\def\shorttitle#1{\def\theshorttitle{#1}}

\let\\\par
\let\thevolumenumber\relax\let\thepapernumber\relax
\let\thevolumeyear\relax\let\thesamplenumber\relax\let\startpage\relax
\let\finishpage\relax\let\publishdate\relax\let\receiveddate\relax
\let\reviseddate\relax\let\accepteddate\relax\let\theasciititle\relax
\let\thecovertitle\relax\let\theasciiauthors\relax\let\theasciiaddress\relax
\let\theasciiabstract\relax\let\theasciikeywords\relax
\let\theasciiemail\relax\let\theshortauthors\relax\let\theshorttitle\relax

\long\def\maketitlep{   

\count0=\startpage

\gt\hfill      
\beginpicture
\setcoordinatesystem units <0.33truein, 0.33truein> point at 2.2 0.9
\setplotsymbol ({$\cal G$})
\plotsymbolspacing=9truept
\circulararc 315 degrees from 0 1 center at 0 0
\setplotsymbol ({$\cal T$})
\circulararc 315 degrees from 1 -1 center at 1 0
\endpicture
\break
{\small\ifx\thesamplenumber\relax 
Volume \else Sample
\fi\thevolumenumber\ (\thevolumeyear)
\startpage--\finishpage\nl
Published: \publishdate}
\vglue 0.5truein plus 0.4fil minus 0.1truein

{\parskip=0pt\leftskip 0pt plus 1fil\def\\{\par\smallskip}{\ifplaintex\large
\else\Large\fi\bf\thetitle}\par\medskip}   

\vglue 0pt plus 0.1fil 

{\parskip=0pt\leftskip 0pt plus 1fil\def\\{\par}{\sc\theauthors}
\par\medskip}

\vglue 0pt plus 0.1fil 

{\small\parskip=0pt\let\newline\\
{\leftskip 0pt plus 1fil\def\\{\par}{\sl\theaddress}\par}
\expandafter\ifx\theemail\relax    
\relax\else\vglue 5pt plus 0.02fil minus 2pt\def\\{\stdspace{\rm 
and}\stdspace} 
\cl{Email:\stdspace\tt\theemail}\fi
\ifx\theurl\relax                  
\relax\else\vglue 5pt plus 0.02fil minus 2pt\def\\{\stdspace{\rm 
and}\stdspace}
\cl{URL:\stdspace\tt\theurl}\fi\par}

\vglue 7pt plus 0.3fil minus 3pt

{\bf Abstract}
\vglue 5pt plus 0.1fil minus 2pt

\theabstract

\vglue 7pt plus 0.3fil minus 3pt

{\bf AMS Classification numbers}\quad Primary:\quad \theprimaryclass

Secondary:\quad \thesecondaryclass

\vglue 5pt plus 0.3fil minus 2pt

{\bf Keywords:}\quad \thekeywords

\vglue 10pt plus 0.5fil minus 5pt

{\small  Proposed: \theproposer\hfill Received: \receiveddate\nl
Seconded: \theseconders\hfill 
\ifx\reviseddate\relax                         
Accepted: \accepteddate                        
\else
Revised: \reviseddate                          
\fi}
\eject
}       

\let\maketitlepage\maketitlep
\let\maketitle\maketitlepage

\font\phead=cmsl9 scaled 950
\font=cmsl9 scaled 1050
\font\pnum=cmbx10 scaled 913
\font\lnum=cmbx10 
\font\pfoot=cmsl9 scaled 950
\font=cmsl9 scaled 1050
\ifplaintex
\headline{\vbox to 0pt{\vskip -4.5mm\line{\small\phead\ifnum
\count0=\startpage ISSN 1364-0380 (on line)
1465-3060 (printed) \hfill {\pnum\folio}\else\ifodd\count0\def\\{ }%
\ifx\theshorttitle\relax\thetitle\else\theshorttitle\fi\hfill{\pnum\folio}
\else\def\\{ and }{\pnum\folio}\hfill\ifx\theshortauthors\relax\theauthors
\else\theshortauthors\fi\fi\fi}\vss}}
\footline{\vbox to 0pt{\vglue 0mm\line{\small\pfoot\ifnum\count0=\startpage
\copyright\ \gtp\hfill\else
\gt, Volume \thevolumenumber\ (\thevolumeyear)\hfill\fi}\vss
}}
\else
\makeatletter
\def\@oddhead{{\small\lhead\ifnum\count0=\startpage ISSN 1364-0380 (on line)
1465-3060 (printed) \hfill {\lnum\number\count0}\else\ifodd\count0
\def\\{ }\ifx\theshorttitle\relax \thetitle \else\theshorttitle\fi\hfill
{\lnum\number\count0}\else\def\\{ and }{\lnum\number\count0}
\hfill\ifx\theshortauthors\relax 
\theauthors\else\theshortauthors\fi\fi\fi}}\def\@evenhead{\@oddhead}
\def\@oddfoot{\small\lfoot\ifnum\count0=\startpage\copyright\ \gtp\hfill\else
\gt, Volume \thevolumenumber\ (\thevolumeyear)\hfill\fi}
\def\@evenfoot{\@oddfoot}
\makeatother
\fi

\newwrite\gtoutfile
\long\gdef\makeheadfile{  
{\def\\{, }\def\s{ }
\immediate\openout\gtoutfile head.xxx
\immediate\write\gtoutfile{To: math@arxiv.org}
\immediate\write\gtoutfile{Subject: put or rep NNNNN:pppp}
\immediate\write\gtoutfile{--text follows this line--}
\immediate\write\gtoutfile{Proxy-for: \ifx\theasciiauthors\relax
\theauthors\else\theasciiauthors\fi\s<\ifx\theasciiemail\relax\theemail\else\theasciiemail\fi>}
\immediate\write\gtoutfile{\noexpand\\}
\immediate\write\gtoutfile{Authors: \ifx\theasciiauthors\relax
\theauthors\else\theasciiauthors\fi}
\immediate\write\gtoutfile{Title: \ifx\theasciititle\relax
\thetitle\else\theasciititle\fi}
\immediate\write\gtoutfile{Subj-class: GT or SG or MG etc}
\immediate\write\gtoutfile{MSC-class: \theprimaryclass\ifx\thesecondaryclass\relax\else, \thesecondaryclass\fi}
\immediate\write\gtoutfile{Journal-ref: Geom. Topol. \thevolumenumber
(\thevolumeyear) \startpage-\finishpage}
\immediate\write\gtoutfile{Comments: Published by Geometry and Topology at}
\immediate\write\gtoutfile{\s\s http://www.maths.warwick.ac.uk/gt/GTVol\thevolumenumber/paper\thepapernumber.abs.html}
\immediate\write\gtoutfile{\noexpand\\}
\immediate\write\gtoutfile{}
\ifx\theasciiabstract\relax
\immediate\write\gtoutfile{\theabstract}\else
\immediate\write\gtoutfile{\theasciiabstract}\fi
\immediate\write\gtoutfile{}
\immediate\write\gtoutfile{\noexpand\\}
\immediate\write\gtoutfile{}
\immediate\closeout\gtoutfile}}  

\def\maketitlepage{\maketitlep\makeheadfile}
\let\maketitle\maketitlepage

\def\ifplaintex{\expandafter\ifx\csname documentclass\endcsname\relax}

\ifplaintex 
\else
\fi

\expandafter\ifx\csname beginpicture\endcsname\relax
\expandafter\ifx\csname documentclass\endcsname\relax
\input pictex \else
\input prepictex \input pictex \input postpictex \fi\fi

\def\gt{{\mathsurround=0pt\it $\cal G\mskip-2mu$eometry \&\ 
$\cal T\!\!$opology}}        

\def\gtp{{\mathsurround=0pt\it $\cal G\mskip-2mu$eometry \&\ 
$\cal T\!\!$opology $\cal P\!$ublications}}  

\def\volumenumber#1{\def\thevolumenumber{#1}}
\def\papernumber#1{\def\thepapernumber{#1}}
\def\volumeyear#1{\def\thevolumeyear{#1}}

\def\pagenumbers#1#2{\def\startpage{#1}\def\finishpage{#2}}
\def\published#1{\def\publishdate{#1}}
\def\proposed#1{\def\theproposer{#1}}
\def\seconded#1{\def\theseconders{#1}}
\def\received#1{\def\receiveddate{#1}}

\def\accepted#1{\def\accepteddate{#1}}
\def\asciititle#1{\def\theasciititle{#1}}
\def\covertitle#1{\def\thecovertitle{#1}}

\def\asciiaddress#1{\def\theasciiaddress{#1}}

\long\def\asciiabstract#1{\long\def\theasciiabstract{#1}}
\def\asciikeywords#1{\def\theasciikeywords{#1}}

\def\shorttitle#1{\def\theshorttitle{#1}}

\let\\\par
\let\thevolumenumber\relax\let\thepapernumber\relax
\let\thevolumeyear\relax\let\thesamplenumber\relax\let\startpage\relax
\let\finishpage\relax\let\publishdate\relax\let\receiveddate\relax
\let\reviseddate\relax\let\accepteddate\relax\let\theasciititle\relax
\let\thecovertitle\relax\let\theasciiauthors\relax\let\theasciiaddress\relax
\let\theasciiabstract\relax\let\theasciikeywords\relax
\let\theasciiemail\relax\let\theshortauthors\relax\let\theshorttitle\relax

\long\def\maketitlep{   

\count0=\startpage

\gt\hfill      
\beginpicture
\setcoordinatesystem units <0.33truein, 0.33truein> point at 2.2 0.9
\setplotsymbol ({$\cal G$})
\plotsymbolspacing=9truept
\circulararc 315 degrees from 0 1 center at 0 0
\setplotsymbol ({$\cal T$})
\circulararc 315 degrees from 1 -1 center at 1 0
\endpicture
\break
{\small\ifx\thesamplenumber\relax 
Volume \else Sample
\fi\thevolumenumber\ (\thevolumeyear)
\startpage--\finishpage\nl
Published: \publishdate}
\vglue 0.5truein plus 0.4fil minus 0.1truein

{\parskip=0pt\leftskip 0pt plus 1fil\def\\{\par\smallskip}{\ifplaintex\large
\else\Large\fi\bf\thetitle}\par\medskip}   

\vglue 0pt plus 0.1fil 

{\parskip=0pt\leftskip 0pt plus 1fil\def\\{\par}{\sc\theauthors}
\par\medskip}

\vglue 0pt plus 0.1fil 

{\small\parskip=0pt\let\newline\\
{\leftskip 0pt plus 1fil\def\\{\par}{\sl\theaddress}\par}
\expandafter\ifx\theemail\relax    
\relax\else\vglue 5pt plus 0.02fil minus 2pt\def\\{\stdspace{\rm 
and}\stdspace} 
\cl{Email:\stdspace\tt\theemail}\fi
\ifx\theurl\relax                  
\relax\else\vglue 5pt plus 0.02fil minus 2pt\def\\{\stdspace{\rm 
and}\stdspace}
\cl{URL:\stdspace\tt\theurl}\fi\par}

\vglue 7pt plus 0.3fil minus 3pt

{\bf Abstract}
\vglue 5pt plus 0.1fil minus 2pt

\theabstract

\vglue 7pt plus 0.3fil minus 3pt

{\bf AMS Classification numbers}\quad Primary:\quad \theprimaryclass

Secondary:\quad \thesecondaryclass

\vglue 5pt plus 0.3fil minus 2pt

{\bf Keywords:}\quad \thekeywords

\vglue 10pt plus 0.5fil minus 5pt

{\small  Proposed: \theproposer\hfill Received: \receiveddate\nl
Seconded: \theseconders\hfill 
\ifx\reviseddate\relax                         
Accepted: \accepteddate                        
\else
Revised: \reviseddate                          
\fi}
\eject
}       

\let\maketitlepage\maketitlep
\let\maketitle\maketitlepage

\font\phead=cmsl9 scaled 950
\font=cmsl9 scaled 1050
\font\pnum=cmbx10 scaled 913
\font\lnum=cmbx10 
\font\pfoot=cmsl9 scaled 950
\font=cmsl9 scaled 1050
\ifplaintex
\headline{\vbox to 0pt{\vskip -4.5mm\line{\small\phead\ifnum
\count0=\startpage ISSN 1364-0380 (on line)
1465-3060 (printed) \hfill {\pnum\folio}\else\ifodd\count0\def\\{ }%
\ifx\theshorttitle\relax\thetitle\else\theshorttitle\fi\hfill{\pnum\folio}
\else\def\\{ and }{\pnum\folio}\hfill\ifx\theshortauthors\relax\theauthors
\else\theshortauthors\fi\fi\fi}\vss}}
\footline{\vbox to 0pt{\vglue 0mm\line{\small\pfoot\ifnum\count0=\startpage
\copyright\ \gtp\hfill\else
\gt, Volume \thevolumenumber\ (\thevolumeyear)\hfill\fi}\vss
}}
\else
\makeatletter
\def\@oddhead{{\small\lhead\ifnum\count0=\startpage ISSN 1364-0380 (on line)
1465-3060 (printed) \hfill {\lnum\number\count0}\else\ifodd\count0
\def\\{ }\ifx\theshorttitle\relax \thetitle \else\theshorttitle\fi\hfill
{\lnum\number\count0}\else\def\\{ and }{\lnum\number\count0}
\hfill\ifx\theshortauthors\relax 
\theauthors\else\theshortauthors\fi\fi\fi}}\def\@evenhead{\@oddhead}
\def\@oddfoot{\small\lfoot\ifnum\count0=\startpage\copyright\ \gtp\hfill\else
\gt, Volume \thevolumenumber\ (\thevolumeyear)\hfill\fi}
\def\@evenfoot{\@oddfoot}
\makeatother
\fi

\newwrite\gtoutfile
\long\gdef\makeheadfile{  
{\def\\{, }\def\s{ }
\immediate\openout\gtoutfile head.xxx
\immediate\write\gtoutfile{To: math@arxiv.org}
\immediate\write\gtoutfile{Subject: put or rep NNNNN:pppp}
\immediate\write\gtoutfile{--text follows this line--}
\immediate\write\gtoutfile{Proxy-for: \ifx\theasciiauthors\relax
\theauthors\else\theasciiauthors\fi\s<\ifx\theasciiemail\relax\theemail\else\theasciiemail\fi>}
\immediate\write\gtoutfile{\noexpand\\}
\immediate\write\gtoutfile{Authors: \ifx\theasciiauthors\relax
\theauthors\else\theasciiauthors\fi}
\immediate\write\gtoutfile{Title: \ifx\theasciititle\relax
\thetitle\else\theasciititle\fi}
\immediate\write\gtoutfile{Subj-class: GT or SG or MG etc}
\immediate\write\gtoutfile{MSC-class: \theprimaryclass\ifx\thesecondaryclass\relax\else, \thesecondaryclass\fi}
\immediate\write\gtoutfile{Journal-ref: Geom. Topol. \thevolumenumber
(\thevolumeyear) \startpage-\finishpage}
\immediate\write\gtoutfile{Comments: Published by Geometry and Topology at}
\immediate\write\gtoutfile{\s\s http://www.maths.warwick.ac.uk/gt/GTVol\thevolumenumber/paper\thepapernumber.abs.html}
\immediate\write\gtoutfile{\noexpand\\}
\immediate\write\gtoutfile{}
\ifx\theasciiabstract\relax
\immediate\write\gtoutfile{\theabstract}\else
\immediate\write\gtoutfile{\theasciiabstract}\fi
\immediate\write\gtoutfile{}
\immediate\write\gtoutfile{\noexpand\\}
\immediate\write\gtoutfile{}
\immediate\closeout\gtoutfile}}  

\def\maketitlepage{\maketitlep\makeheadfile}
\let\maketitle\maketitlepage

\volumenumber{4}\papernumber{14}\volumeyear{2000}
\pagenumbers{407}{430}
\proposed{Robion Kirby}
\seconded{Wolfgang Metzler, Cameron Gordon}
\received{22 May 2000}
\accepted{3 November 2000}
\published{10 November 2000}

\usepackage{epsf,amssymb} 

\newtheorem*{theorem} {Theorem} 
\newtheorem*{corollary} {Corollary} 
\newtheorem*{lemma} {Lemma} 
\def\S{Section }

\newcount\sectionnumber            
\newcount\resultnumber             
\def\Section#1{\global\advance\sectionnumber by 1
\vskip-\lastskip\penalty-800\vskip 20pt plus10pt minus5pt 
{\Large\bf\number\sectionnumber\quad#1}         
\vskip 8pt plus4pt minus4pt
\nobreak\resultnumber=1}      
\def\proc#1{
\vskip-\lastskip\ppar
{\bf\number\sectionnumber.\number\resultnumber
\qua#1\qua}\global\advance\resultnumber by 1\ignorespaces}

\newcommand{\figureone}{
$$\begin{picture}(200,175)  \small
    \put(-65,0)       {\epsfbox{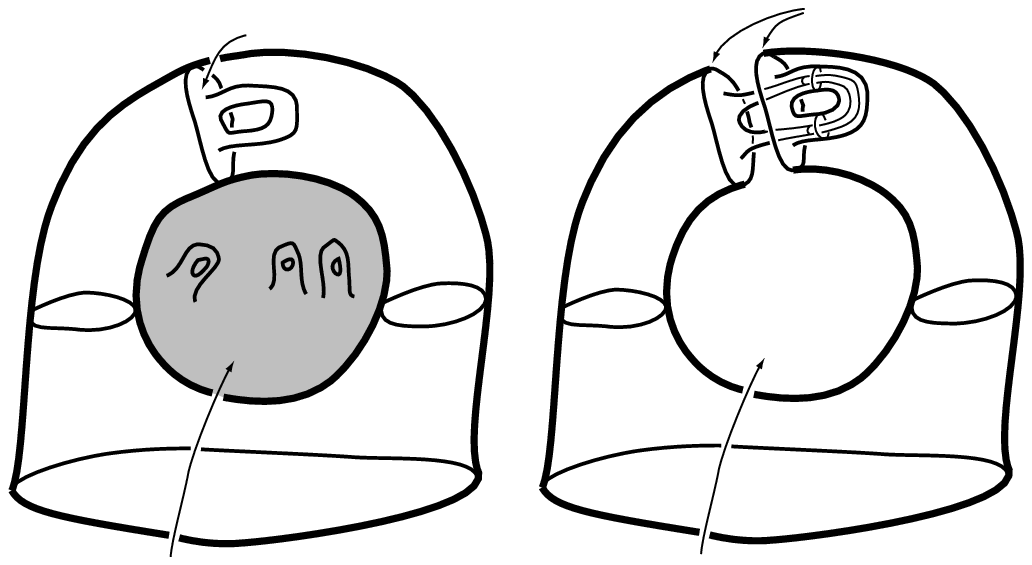}}
    \put(10,156) {subgrope}
    \put(-20,-6) {dual subgrope}
    \put(174, 157) {cut, glue parallels}
    \put(130,-6) {delete}
    \end{picture}$$
}

\newcommand{\figuretwo}{
$$\begin{picture}(200,155)  \small
    \put(-90,-5)       {\epsfbox{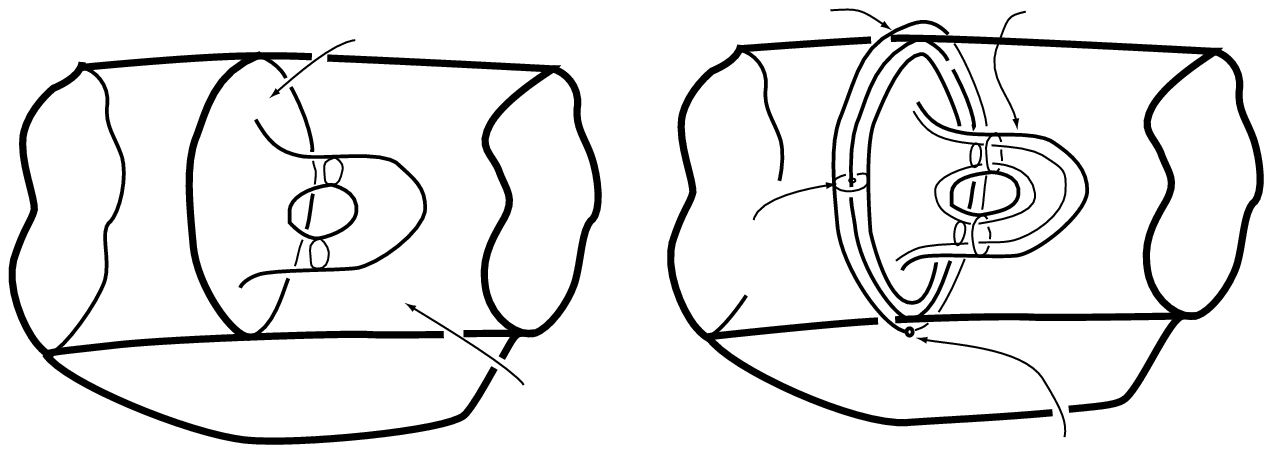}}
    \put(-50,121)   {$H$, subgrope dual to $A$}
    \put(-30,13)   {$A$}
    \put(66,12)   {$B$}
    \put(95,131)   {torus around} 
    \put(135,121)   {$\partial H$}
    \put(115,58)   {normal}
    \put(130,48)   {disk}
    \put(130,-6)   {intersection of torus with $A$}
    \put(175,133)   {parallel of $H$, attached}
    \put(215,123)   {to torus}
    \end{picture}$$
}

\newcommand{\figurethree}{
$$\begin{picture}(200,175) \small
    \put(-60,0)       {\epsfbox{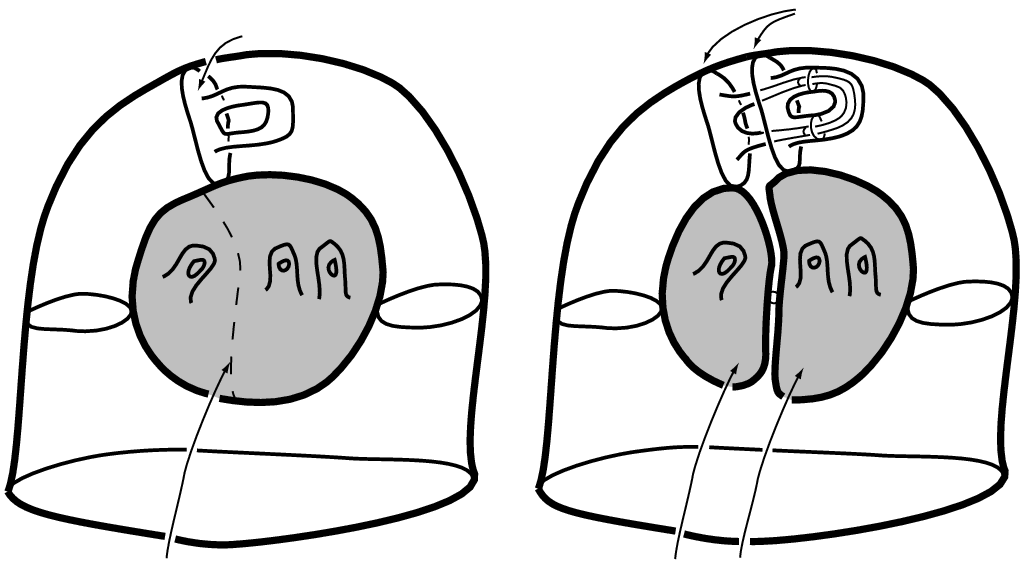}}
    \put(14,155) {$H$, dual to $A$}
    \put(-22,-6) {arc on $A$}
    \put(175, 159) {parallel copies of $H$}
    \put(130,-7) {$A$, split}
    \end{picture}$$
}

\newcommand{\figurefour}{
$$\begin{picture}(200,185) \small
    \put(-80,0)       {\epsfbox{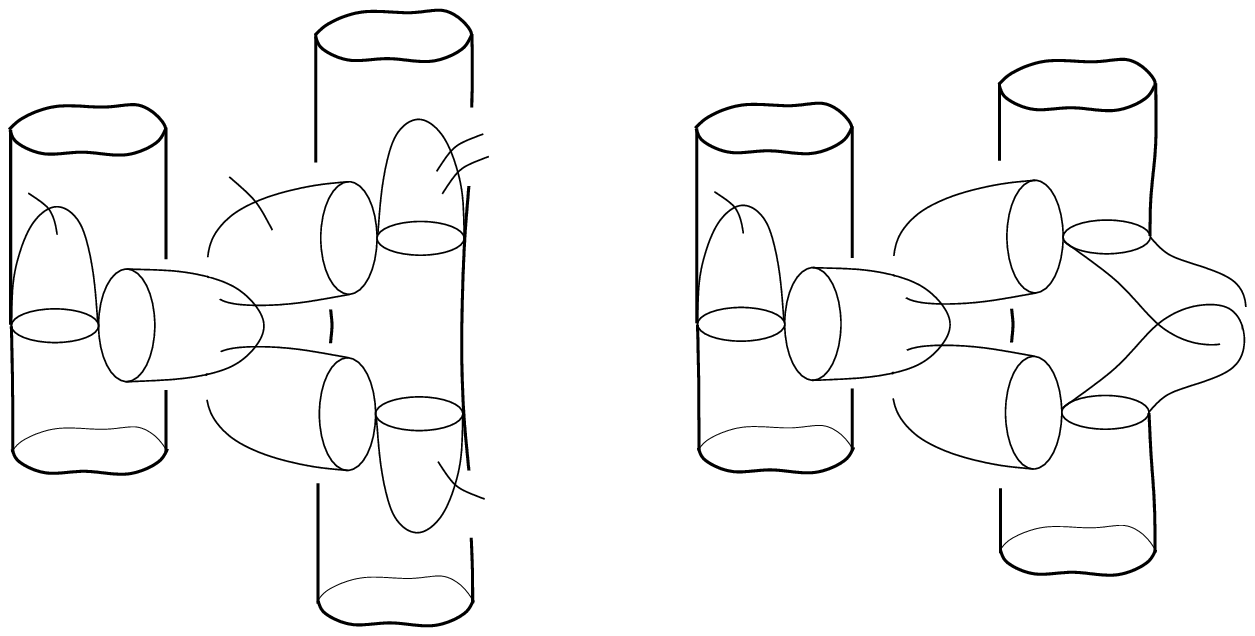}}
    \put(-60,32)      {$B$}
    \put(-14,132)    {$x$}
    \put(-7,57)        {$y$}
    \put(0,100)        {$y$}
    \put(-36,75)      {$x$}
    \put(33,40)      {$z$}
    \put(33,131)    {$z$}
    \put(61,36)        {$w$}
    \put(62,135)     {$v$}
    \put(59,144)     {$v$}
    \put(191,57)     {$y$}
    \put(198,100)    {$y$}
    \put(260,61)     {$z$}
    \put(260,111)    {$z$}
    \end{picture}$$
}

\begin{document}

\title{Subexponential groups in ${\mathbf 4}$--manifold topology}
\shorttitle{Subexponential groups in 4-manifold topology}
\covertitle{Subexponential groups in 4--manifold topology}
\asciititle{Subexponential groups in 4-manifold topology}

\authors{Vyacheslav S Krushkal\\Frank Quinn} 

\address{Department of Mathematics,Yale University\\New Haven, CT 06520-8283,
USA}
\secondaddress{Department of Mathematics,Virginia Tech\\Blacksburg, 
VA 24061-0123, USA}
\asciiaddress{Department of Mathematics,Yale University\\New Haven, 
CT 06520-8283, USA\\Department of Mathematics,Virginia Tech\\Blacksburg, 
VA 24061-0123, USA}
\email{krushkal@math.yale.edu, quinn@math.vt.edu}

\begin{abstract} 

We present a new, more elementary proof of the Freedman--Teichner
result that the geometric classification techniques (surgery,
s--cobordism, and pseudoisotopy) hold for topological 4--manifolds with
groups of subexponential growth.  In an appendix Freedman and Teichner
give a correction to their original proof, and reformulate the growth
estimates in terms of coarse geometry.

\end{abstract}

\asciiabstract{We present a new, more elementary proof of the
Freedman-Teichner result that the geometric classification techniques
(surgery, s-cobordism, and pseudoisotopy) hold for topological
4-manifolds with groups of subexponential growth.  In an appendix
Freedman and Teichner give a correction to their original proof, and
reformulate the growth estimates in terms of coarse geometry.}

\primaryclass{57N13}
\secondaryclass{57N37, 57N70, 57R65}

\keywords{$4$--manifolds, groups of subexponential growth, gropes}
\asciikeywords{4-manifolds, groups of subexponential growth, gropes}

\maketitlepage

The disk embedding theorem for 4--manifolds with ``good'' 
fundamental group is the key ingredient 
of the 
classification theory: it is used in the proof of the 4--dimensional 
surgery theorem, and 
the 5--dimensional s--cobordism theorem and pseudoisotopy 
theorems. The 
homotopy hypotheses 
of the theorem always allow one to find a 2--stage immersed 
capped grope. If one can find   
such a 
grope so that loops in the image are nullhomotopic in the 
ambient manifold, then 
Freedman's theorem \cite{F, FQ} shows 
there is a topologically flat embedded disk. 
The current focus, therefore, is on obtaining this $\pi_1$--nullity 
condition. Freedman \cite{F1} 
showed this is possible if the fundamental group of the manifold is 
poly-(finite or cyclic). 
This was extended to groups of polynomial growth in \cite{S}. 
The current best result is 
for 
groups of subexponential growth. 

The disk theorem for subexponential groups was stated by 
Freedman and Teichner \cite{FT}. However the ``key point'' of
\cite{FT} page 521, line 17, is incorrect.
In the Appendix  Freedman and Teichner show how to
modify their construction to correct this. 
The present paper sidesteps the issue by using
a different and more elementary construction 
developed by the first author (see \cite{K}). It displays 
particularly
clearly how the proof fails in 
the general (exponential growth) case, and suggests that 
an infinite construction may
be necessary 
to make further progress.  

The following result is the input needed for the  disk 
embedding theorem 
for manifolds with subexponential
fundamental groups. For a full statement of the disk 
embedding theorem and applications to surgery and 
s--cobordism,
see~\cite{FT}. For the application 
to pseudoisotopy see~\cite{Q}.

\begin{theorem} \sl
Suppose $G\to M^{4}$ is a properly immersed (capped) 
grope of height $\geq 2$, 
and $\rho\co \pi_{1}M \to \pi$ is a homomorphism with 
$\pi$ of subexponential growth.  
Then the total contraction of $G$ is regularly 
homotopic rel boundary 
to an immersion whose double point loops have 
trivial image in $\pi$.
\end{theorem} 

This slightly extends the usual immersion-improvement 
formulation in that we do not 
require the total contraction to be a disk, and the output 
immersion is regularly homotopic 
to the input. Neither extension has new consequences, 
but they come for free in the 
proof and they simplify applications. Capped gropes, 
contractions, and subexponential 
growth are all reviewed in the text.

In rough outline the proof goes as follows: The images 
of the double point loops 
of $G$ give a finite subset of $\pi$. Subexponential 
growth implies that in a large 
collection of words of fixed length in the finite subset, 
somewhere there is a subword 
whose product is trivial. We organize the data so this 
subword can be realized 
geometrically as double point loops, by pushing 
intersections around in the grope. This 
eliminates a branch in the grope in exchange for 
$\pi$--trivial self-intersections of 
the base surface. Iterating this eliminates all 
branches (ie, gives a contraction) with 
$\pi$--trivial self-intersections. The key 
technique is splitting to dyadic branches.

{\bf Acknowledgements}\qua VK was supported at the Institute 
for Advanced Study by NSF grant 
DMS 97-29992. FQ is partially supported by NSF grant 
DMS 97-05168 and Harvard University.

\Section{Definitions} 

We briefly review the definition of gropes in order to 
fix terminology for caps, 
duals, immersions, contractions, etc. 

\proc{Gropes}
Begin with a surface $S$.  
A {\it model grope\/} built from $S$, or an 
{\it $S$--like grope\/} is a 2--complex in $S\times I$ 
obtained by repeatedly replacing embedded disks by 
punctured tori with disks attached, see 
\cite[\S 2.1]{FQ}.  The attached disks are the ``caps'' 
for the torus; there are two of them, they 
intersect in a point, and each is referred to as the 
{\it dual\/} of the other. 

\begin{enumerate}
\item  The {\it caps\/} of the final grope are caps 
introduced at some stage, and which have not been 
modified later. 
\item A {\it subgrope\/} is a cap that has been modified,
 together with these modifications.  
More generally a subgrope is a disjoint union of these. 
Subgropes are disk-like gropes, or more 
generally (union of disk)-like gropes.
\item {\it Dual\/} subgropes are subgropes obtained 
from caps 
that were dual at some stage of the construction. 
They are attached to the same lower surface, 
and the
attaching circles intersect in a single point.
\item The {\it base\/} of the grope is the 
surface obtained by modifying $S$.
\item A {\it branch\/} is a dual pair of subgropes 
attached to 
the base.
\item The grope has {\it height\/} $\geq k$ if any path 
from a cap to the base passes 
through at least $k$ boundary curves of subgropes. 
In other words there are at least $k$ levels 
of surfaces (counting the base) below each cap. 
Note each branch has height $\geq k-1$.
\end{enumerate}

\proc{Contractions}
Suppose $H_{+}, H_{-}$ are dual subgropes 
in a grope. The {\it contraction away from\/} 
$H_{-}$, or equivalently the {\it contraction across\/}
$H_{+}$ is defined as follows: discard
the  subgrope $H_{-}$, cut open the surface to 
which $H_{+}$
is attached  along the attaching circle, and glue in two
parallel copies of $H_{+}$. See Figure~1. This
process is sometimes referred to as surgery. The 
result is an $S$--like grope in a canonical way, ie, 
there is a canonical way to see 
this as obtained from $S$ in $S \times I$.

\begin{figure}[ht]
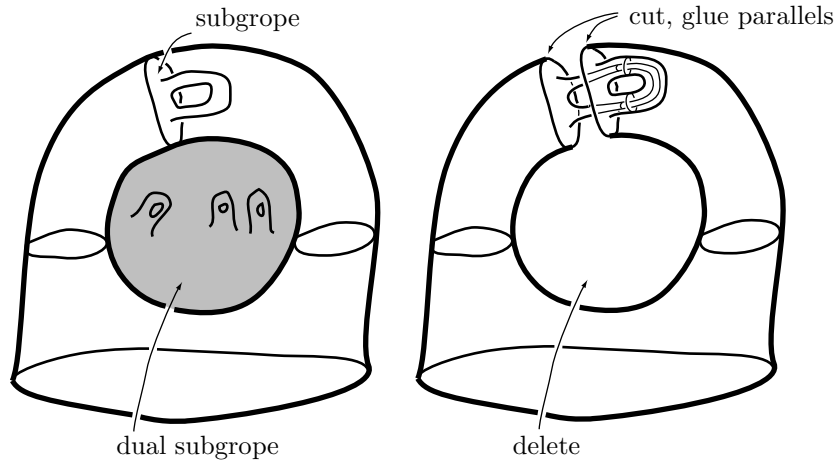
 
\figureone
\caption{Contraction}
\end{figure}

There is also a ``symmetric contraction'' which cuts 
the surface along the attaching curves 
of both subgropes and glues in parallel copies of both, 
together with a small square near 
the intersection point of the curves. See \cite[\S 2.3]{FQ}. 
This operation won't be used here.

\begin{lemma} \sl 
Suppose $H_{+}, H_{-}$ are dual caps in an $S$--like 
grope. Then the contractions 
across $H_{+}$, across $H_{-}$, and the symmetric 
contraction, give isotopic gropes.
\end{lemma} 

\proof
Here ``isotopic'' means ambient isotopic in
$S\times I$. Since $H_{+}, H_{-}$ are both caps, all 
these contractions undo one disk-punctured 
torus replacement. The isotopy is clear in pictures; 
see \cite[\S 2.3]{FQ}. 
\endproof

A {\it total contraction\/} of a grope is obtained by 
repeatedly contracting (in any order) until no caps
are left. Iterating the lemma shows this always 
returns the original surface:

\begin{corollary} \sl 
All total contractions of $S$--like gropes are surfaces 
isotopic (rel boundary, in $S\times I$) to $S$.\endproof
\end{corollary} 

When a grope is immersed the different contractions 
give regularly homotopic
immersions of $S$. Their usefulness is that they may 
differ greatly in their intersection patterns.

\proc{Immersions}
A {\it proper immersion\/} of a grope $G\to M^{4}$ is an 
immersion of a regular neighborhood of the spine 
in $S\times I$ satisfying

\begin{enumerate}
\item boundary goes to boundary, so the image 
intersects  $\partial M$
in $\partial S\times I$;
\item all intersections come from transverse 
intersections among caps; and
\item (if the base surface is noncompact) the 
immersion is proper in the topological sense.
\end{enumerate}

Caps have $I$--bundle neighborhoods in
$S\times I$. Condition (2) means these neighborhoods 
intersect in squares determined by
transversality and the bundle structures, and that 
there are no other intersections. 

In fact a regular neighborhood of the grope spine  is 
isotopic to $S\times I$ itself. The 
definition is given in terms of regular neighborhoods so the 
standard $D^{2}$ bundle 
neighborhoods 
of caps will be easier to see. A consequence of this 
neighborhood uniqueness is that the 
original core copy of $S$ is isotopic into any regular 
neighborhood of the grope spine. 
Composing this with a grope immersion gives an 
immersion of $S$, well-defined up to 
regular homotopy. As remarked above, total 
contractions give explicit descriptions of such 
immersions.

\proc{Transverse spheres and gropes}
A {\it transverse sphere\/} for a surface in a 
4--manifold is a framed immersed sphere that intersects 
the surface in a single point, see \cite[\S 1.9]{FQ}. 
Similarly a transverse grope is an immersed
sphere-like grope that intersects the surface in 
one point, and this point is in the base of the grope.
Note that totally contracting  a transverse grope 
gives a transverse sphere. 

In a grope, every 
surface component not part of the base has a 
standard transverse grope. This is constructed in 
$S\times D^{2}$ ($S$ is the original surface, 
or the total contraction), so there is a copy inside 
any neighborhood of the spine. This was a key 
part of Freedman's original constructions of 
convergent gropes. The {\it tight\/}  transverse 
grope for a surface component $A$ is obtained as
follows. Let $B$ be the next surface down, so 
$A$ is attached to $B$ along a circle. Let $H$ be 
the subgrope dual to $A$. The base of the 
transverse grope is the torus obtained as the normal 
$S^{1}$ bundle to $B$, restricted to a parallel 
of the attaching curve of $H$. To this we add a 
parallel copy of $H$ itself, and a $D^{2}$ fiber 
of the normal disk bundle of $B$. The result is a 
sphere-like grope, see Figure 2. Connected sum 
with this grope gives the ``lollipop'' move of 
\cite{FT}.

\begin{figure}[ht]
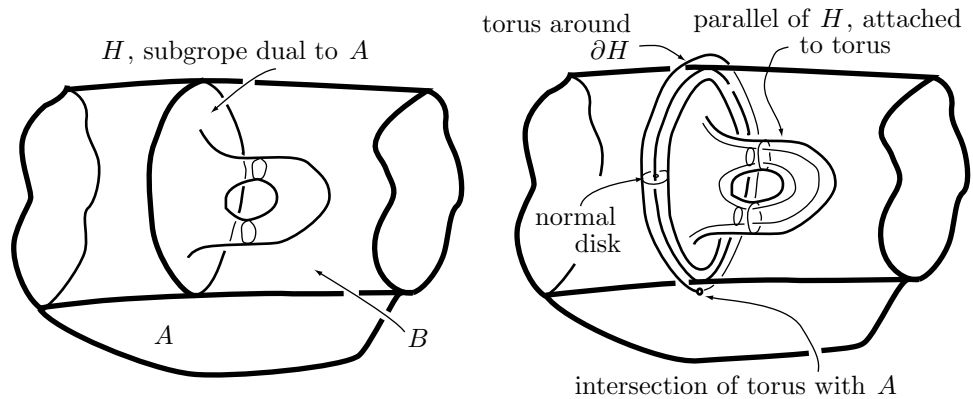
 
\figuretwo
\caption{Tight transverse grope}
\end{figure}

The tight transverse grope has two points of 
intersection with the grope: 
one is the (desirable) intersection between the 
base torus and $A$. The other is the intersection 
of the normal-disk cap with $B$. Usually we 
avoid this one by contracting away from this cap. 
(See 1.2. Recall that this discards the normal-disk 
cap, cuts the torus along 
the boundary of the parallel of $H$, and fills in 
with two copies of $H$.) This contraction is the 
``loose'' or standard transverse grope. It is 
described in \cite[\S 2.6]{FQ}, and sums with it give 
the ``double lollipop'' move of \cite{FT}. 

Finally, transverse spheres for $A$ are obtained 
by totally contracting the 
loose transverse grope. 

We repeat here the warning of \cite[\S 2.6]{FQ}. 
If $H_{+}$, $H_{-}$ 
are dual subgropes then the construction gives 
transverse gropes for the bottom-level surfaces 
in both. However these two gropes intersect each 
other, which in practice  means they cannot 
be used simultaneously. (This is a mistake in \cite{FT}.) 
This problem can be partially avoided 
by using them sequentially, see \cite[\S 2.7]{FQ}. 

\begin{lemma} \sl The spheres obtained by 
totally contracting loose
transverse gropes bound embedded 3--disks in 
$S\times D^{2}$. Thus 
connected sum with these spheres changes a 
surface by regular homotopy.
\end{lemma} 

\proof
To see this begin with the tight grope and totally 
contract the copy of $H$. This gives a torus with two 
caps. The torus bounds a solid torus (the $D^{2}$ 
bundle over the attaching circle of $H$) disjoint 
from the $H$ cap. Attach to this solid torus a 
$D^{2}\times I$ thickening of the $H$ cap. The  
result is a 3--disk whose boundary is the transverse sphere.
\endproof

\Section{Dyadic branches and splitting}

A disk-like grope is {\it dyadic\/} if all component 
surfaces are either disks 
or punctured tori. This means each non-cap surface 
has exactly one pair of dual subgropes attached 
to it. A grope has {\it dyadic branches\/} if all of the 
subgropes above the base level are dyadic. 
The benefits of dyadic branches are especially simple 
contractions and simpler tracking of their
interactions. In 2.3 we describe the
splitting construction, which converts any grope into 
one with
dyadic branches. Later in the section it is used further 
to set up a (dyadic) situation in which 
product information in $\pi$ can be exploited.

\proc{Contraction of dyadic gropes}
When a grope is contracted caps are discarded. In a total 
contraction most get 
discarded. Dyadic gropes are special in that there are
total contractions in which {\it all but one\/} 
of the caps are discarded. Since intersections occur 
among caps, discarding them simplifies 
intersection data. Discarding all but one gives the 
greatest possible simplification of the data. 

To be explicit, suppose $C$ is a cap in a dyadic branch. 
The {\it total contraction across $C$} is 
defined as follows: $C$ lies in one of the dual 
subgropes in the branch; form the contraction 
across this subgrope. The result has two dyadic 
branches parallel to the subgrope, so each 
contains a parallel copy of $C$. Repeat this 
construction, contracting each new branch across 
the subgrope containing the copy of $C$. Each 
iteration doubles the number of copies of $C$, 
and reduces their height above the base by 1. 
If the original height of $C$ is $k$, then  total 
contraction eliminates the branch and subsitutes $2^{k}$  
parallel copies of $C$ in the base 
surface.

\proc{Pushing down and back up}
Suppose $C_{1}$, $C_{2}$ are 
distinct caps in a dyadic branch in a grope. 
We construct a transverse sphere for $C_{1}$ that 
contains parallel copies of $C_{2}$, and otherwise 
lies in a small neighborhood of the body of 
the grope. Suppose a surface $W$ intersects 
$C_1$. Then adding copies of this transverse 
sphere to remove the intersection points will 
be called {\it pushing $W$ down off $C_{1}$ and 
back up across $C_2$}. Reasons for the terminology are:

\begin{enumerate}\item ``pushing'' because the 
new surface is regularly homotopic to $W$;
\item ``off $C_{1}$'' because intersections with 
$C_{1}$ have been removed; and
\item ``across $C_{2}$'' because the modified 
surface contains 
parallel copies of  $C_{2}$, and will therefore 
intersect anything $C_{2}$ intersects.
\end{enumerate}

We now construct the transverse sphere. Let 
$H_{1}$ and $H_{2}$ be dual subgropes in the branch, 
so that $H_{1}$ contains   $C_{1}$ and
$H_{2}$ contains $C_{2}$. Let $\hat T$ be the  tight 
transverse grope for the bottom surface of $H_{1}$. 
This is a dyadic sphere-like grope with a 
cap parallel  to $C_{2}$; let $T$ be the sphere 
obtained by totally contracting $\hat T$ across 
this cap.   $T$ contains $2^j$ parallels of $C_{2}$, 
where $j$ is the height of $C_{2}$ in $H_2$.
If $H_1=C_1$ then this is the desired sphere. 
Otherwise take a 2--sphere fiber of the normal 
sphere bundle of the attaching circle of $C_{1}$. 
This 2--sphere intersects $C_{1}$ in one point 
and intersects
the surface to which it is attached in two points. 
Push these latter points down 
in the grope to the bottom surface of $H_{1}$. If 
$C_{1}$ has height $k$ in $H_{1}$ then there 
will be $2^{k}$ intersection points there. Remove 
all these by connected sum with parallels of 
the sphere $T$. This gives the desired transverse 
sphere for $C_{1}$. Note that all together this 
sphere has $2^{j+k}$ copies of~$C_{2}$.

This transverse sphere bounds a 3--disk in the model. 
To see this, first note that $T$ does by Lemma 1.4. 
The fiber of a 2--sphere bundle bounds the 
fiber of the associated 3--disk bundle, and the isotopy 
pushing the boundary down in $H_{1}$ 
pushes the 3--disk too. Finally the boundary sum of a 
bunch of 3--disks gives  a 3--disk whose 
boundary is the connected sum of spheres. It is 
easy to see that the 3--disk obtained this way 
is embedded.

Since the sphere bounds an embedded 3--disk, 
connected sum with the sphere 
changes surfaces by regular homotopy.

\proc{Splitting gropes}
The splitting operation splits a surface into two 
pieces, at the cost of doubling the dual subgrope. 
It can be used to decompose branches into dyadic 
branches, and can separate intersection points
distinguished by properties unaffected by the dual doubling. 

Suppose $A$ is a component surface 
of a grope, not part of the base. Let $B$ be the 
surface it is attached to, and $H$  the dual subgrope. 
Now suppose $\alpha$ is an embedded arc on $A$, 
with endpoints on the boundary, and disjoint from
attaching circles of higher stages. In the 3--dimensional 
model, sum $B$ with itself by a tube about
$\alpha$ (the normal $S^{1}$ bundle), and discard 
the part of $A$ that lies  inside the tube. This splits
$A$ into two components. $H$ is a dual for one 
component; obtain a dual for the other by taking a
parallel copy of $H$. See Figure~3.

\begin{figure}[ht]
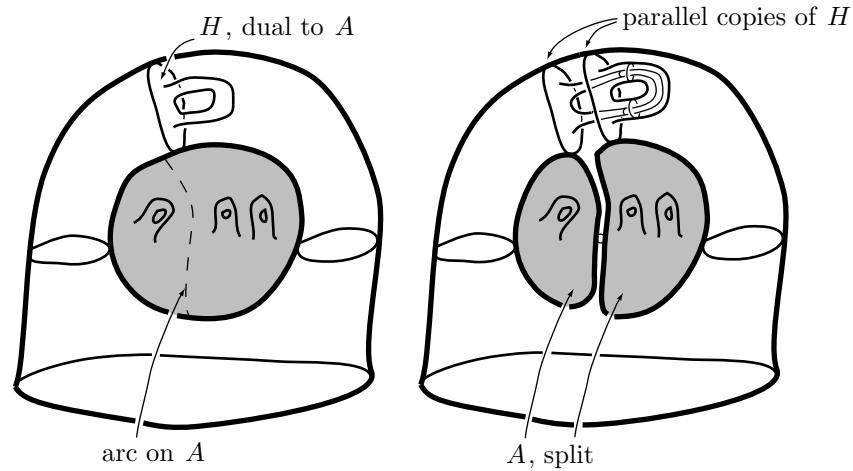
 
\figurethree
\caption{Splitting}
\end{figure}

\begin{lemma} \sl 
Any grope can be transformed by
iterated splitting to one with dyadic branches. If no 
caps are split, and all caps in the result are 
parallel copies of the original caps, then the result is 
well-defined up to isotopy.
\end{lemma} 

\proof
The splitting is done by induction downward from the 
caps. Suppose all subgropes of height $\leq k$ 
that do not include a component of the base are 
dyadic. To start the induction note that this is 
true for $k=0$. For the induction step choose 
splitting arcs on base surfaces of subgropes of height 
$k+1$, provided this surface is not part of the global 
grope base. Choose the arcs so that each
component of the complement has genus 1. Splitting 
one surface may double another, so we 
choose an order to ensure the process terminates. 
First split all surfaces whose dual subgrope 
has height $\leq k$. These duals are already dyadic, 
so doubling does not introduce nondyadic
subgropes. Second, split one in each dual pair where 
both subgropes have height $k+1$. Nondyadic 
subgropes are introduced by doubling the dual. 
However these all now have dyadic duals. Thus 
splitting them does no further harm. After that we 
proceed by induction in $m$, splitting surfaces 
whose duals have height $k+m$. Previous steps ensure 
the height--$(k+1)$ subgropes of the 
duals are already dyadic, so the doubling of the 
dual does not introduce new nondyadic subgropes of
that height. This proves the existence of dyadic 
splittings.\endproof

This process can be seen as an
expansion of a product of sums. Suppose $H_{1}$ 
and $H_{2}$ are dual subgropes, and each has 
a base surface of genus 2. This means each has
two branches, say $H_{i,1}$ and $H_{i,2}$. Think 
of a grope as a sum of its branches, so $H_{1}=
H_{1,1}+H_{1,2}$. Think of a branch as a product 
of the two dual subgropes, so the branch formed 
by the $H_{i}$ becomes 
$$H_{1}*H_{2}=(H_{1,1}+H_{1,2})*(H_{2,1}+H_{2,2}).$$
Splitting converts this into four branches:
$$H_{1,1}*H_{2,1}+H_{1,2}*H_{2,1}+H_{1,1}*H_{2,2}+
H_{1,2}*H_{2,2}.$$
This formulation extends 
to represent an entire grope as an iterated composition 
of polynomials in its caps. Splitting
corresponds to rewriting this composition as a sum of 
monomials. We caution, however, that the
``algebra'' of caps is not associative or commutative.

The uniqueness for minimal dyadic splittings 
can easily be proved by formally
following the proof of uniqueness of monomial expansions 
of iterated 
polynomials. We omit this since we have no application 
for it.

\proc{Splitting and labeling caps}
Suppose there are several different ``types'' of 
intersection points on a 
cap. We can separate the types by arcs, then split 
along the arcs to get caps containing only one 
type of point. Other caps get doubled during this 
process, so for it to succeed ``type'' must be 
defined so that new intersections with parallel copies 
of a cap have the same ``type.'' We apply 
this principle to get uniform local patterns of 
intersection invariants.

\rk{DATA}
To facilitate the definition of invariants we now 
require every 
grope branch to have {\it dyadic labels\/}, and  
{\it a path to the basepoint}. (This last is not 
necessary if the grope is modeled on a simply-connected 
surface, ie, 
$D^{2}$--like, or $S^{2}$--like).

\medskip

Dyadic labels are obtained as follows: for each 
dual pair of surfaces, 
label one by 0 and the other by 1. A cap gets a label 
by reading off the sequence of 0 or 1s 
encountered in a path going from the base to the 
cap. In a dyadic branch the labels uniquely 
specify the caps, since at each level there are 
only two ways to go up and these are 
distinguished by the labels.  When branches are 
split the fragments inherit labelings and paths 
to the basepoint in the evident ways, so this 
data is preserved.

The following is the first 
application of splitting to simplify cap types. It will 
be iterated in later constructions.

\begin{lemma} \sl A grope whose branches have 
dyadic labels and paths to the basepoint can 
be split so that it has dyadic branches and each 
cap satisfies:

\begin{enumerate}
\item there are no self-intersections;
\item all caps intersecting the given one have the 
same label; and 
\item the fundamental group classes of the loops 
through the intersection points are the same.
\end{enumerate}

Further, the subset of $\pi_{1}M$ occuring as 
double point loops is the same as that 
of the original grope.
\end{lemma} 

We clarify that condition (2) requires all the intersecting 
caps 
to have the same label, but this may be different from 
the label of the cap being intersected. 
The loops in (3) go from the fixed cap to the others. 
In detail they go along the path from the 
basepoint and up through the branch to the fixed cap, 
through the intersection point to the other 
cap, then back down and back to the basepoint. Note 
that viewing an intersection point from the 
other
cap reverses this path, so gives the inverse element in 
$\pi_1M$.

\proof In each cap
choose arcs to separate intersection points into sets, 
so that self-intersection points are separated,
and all points in a set have the same group element 
and intersect caps with the same label. 
Initially these sets may each contain only one point. 
Now split to separate all these sets, and 
continue splitting to obtain dyadic branches. If 
parallel copies of a cap are used before it is 
split then use parallels of the chosen splitting arcs 
in these copies. An argument similar to the 
proof of Lemma 2.3 shows that this process terminates. 

The doubling involved in splitting 
means that each original intersection point bifurcates 
to many intersections with parallels of 
pieces of the other cap. However by the way 
labels are chosen in splittings all these pieces 
will have the same label as the original. Similarly 
they have the same group element. Therefore 
the end result has the properties specified in the 
lemma.
\endproof

\proc{Intersection types}
We iterate the construction of 2.4 to arrange all 
intersections with a given cap 
to have the same ``type'' in  more elaborate 
senses.~{The} first 
definition formalizes the situation of 2.4, the 
next inductively extends it.

\rk{Definition} \sl 1--types\qua\rm
\begin{enumerate}
\item A {\it $1$--type\/} is a function from a set 
of dyadic labels to pairs $(\alpha, x)$ with
$\alpha\in\pi_{1}M$ and $x$ a dyadic label;
\item a dyadic branch {\it has 1--type\/} $\rho$ if 
\begin{enumerate}
\item the domain of $\rho$ is the set of labels on 
caps in the branch; and
\item  if $y$ is a cap label and $\rho(y)=(\alpha, x)$, 
then all caps intersecting 
$y$ have label $x$, and all intersection points 
have $\pi_{1}$ element $\alpha$. 
\end{enumerate}
\end{enumerate}

The output of Lemma 2.4 is a grope each of whose 
branches has a 1--type.

\rk{Definition}\sl$n$--types\qua\rm
If $n>1$ then 
\begin{enumerate}
\item An  {\it $n$--type\/} is a function from a set of 
dyadic labels to pairs 
$(\alpha, x)$ with $\alpha$ an $(n-1)$--type and 
$x$ a dyadic label;

\item a dyadic branch {\it has $n$--type\/} $\rho$ if 
\begin{enumerate}
\item the domain of $\rho$  is the set of labels on 
caps in the branch; and
\item  if $y$ is a cap label and $\rho(y)=(\alpha, x)$, 
then all caps 
intersecting $y$ have label $x$, and the branches 
containing these caps have 
$(n-1)$--type $\alpha$. 
\end{enumerate}
\end{enumerate}

Branches with 1--types have uniform intersections 
with adjacent branches. 
Bran\-ches with  $n$--types have uniform patterns 
of intersection going out through chains of $n$
branches, see Figure 4.

\begin{figure}[ht]
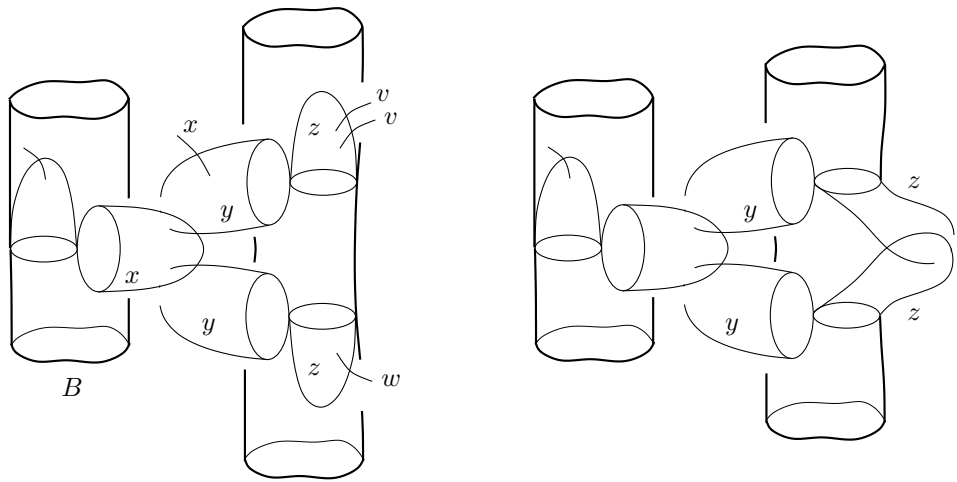
 
\figurefour
\caption{Example: the branch $B$ of a capped surface 
has a $2$--type if 
(1) all ``$x-y$'' intersections determine the same 
element in ${\pi}_1 M$,
(2) the dyadic cap labels $v$ and $w$ coincide, 
and (3) all ``$z-v$'' and ``$z-w$'' intersections 
determine the same element in ${\pi}_1 M$. 
The second figure shows a collision at distance $2$.}
\end{figure}

\eject

\begin{lemma} \sl \hfill
\begin{enumerate}
\item If a branch with an $n$--type is split into 
dyadic branches,
then each of these has the same $n$--type.
\item Any grope can be split to one in which every  
branch has an $n$--type. 
\end{enumerate}
\end{lemma} 

\proof The first statement is 
straightforward. Statement (2) proceeds by induction, 
with Lemma 2.4 starting the induction with
$n=1$.

Suppose that each branch has an $(n-1)$--type. 
Associate to each intersection
point on 
a cap the $(n-1)$--type of the branch intersected. 
Choose arcs in the cap to separate these 
points into sets, so that all the points in each set 
have the same type. Do this for every cap, 
then split along these arcs and continue splitting 
to get dyadic branches. Note that each 
intersection point in the original may bifurcate into 
a large number of intersections with 
fragments of splittings. However since (by (1)) 
all these fragments have the same $(n-1)$--type, 
we end up with all branches intersecting a cap 
having the same $(n-1)$--type.In other words 
the branch has an $n$--type, and the conclusion 
of the lemma
holds.
\endproof

\proc{Collisions}
The construction of 2.4 separates self-intersections and 
makes $\pi_{1}$ classes of
intersections uniform. The construction of 2.5 
extends the uniformity 
out to distance $n$ from each branch. Here we 
similarly extend the self-intersection condition 
out distance $n$, see Figure 4.

\rk{Definition}
A {\it collision at distance\/} $n$ from 
a branch $B$ is a sequence of dyadic labels 
$(x_{i},y_{i})$  and two sequences of branches 
$B_{1,i}$, $B_{2,i}$ for $1\leq i\leq n$ so that

\begin{enumerate} 
\item $B= B_{1,1}=B_{2,1}$
\item for each $i<n$ and $j=1$ or 2, there is an 
intersection between the $x_{i}$ cap of $B_{j,i}$ 
and the $y_{i}$ cap of $B_{j,i+1}$; and
\item the $x_{n}$ cap of $B_{1,n}$ transversally 
intersects the $x_{n}$ cap of
$B_{2,n}$.
\end{enumerate}

We will see that if the branches have $n$--types 
then the only ``unexpected''
intersections that can occur must
satisfy these conditions.

\begin{lemma} \sl Given $n$ and a grope, there is 
a splitting so that there are no collisions at      
distance$\leq n$.
\end{lemma} 

\proof
This proceeds by induction in $n$. A collision at 
distance 1 is just a self-intersection in a cap, so
Lemma 2.4 starts the induction with $n=1$. Suppose 
there are no collisions at distance
$n-1$.
If there is a collision at distance $n$ then there are 
sequences of branches and labels as specified 
in the definition. If $B_{2,1}=B_{2,2}$ then there 
would be a collision at distance
$n-1$, so 
the induction hypothesis implies these are distinct. 
Their intersection points with the cap in $B$ 
must
also be distinct. If we split the cap to separate these 
points we avoid the collision, or more
precisely, we postpone
it to  distance $n+1$. Moreover splitting does not 
introduce new collisions 
at shorter distances. The induction step thus
proceeds by choosing  arcs to separate any points that
lead  to a collision at distance $n$, and splitting along
these.
\endproof

\Section{Groups with subexponential growth}
A group $\pi$ is said to have {\it subexponential growth\/} 
if given any
finite subset $S\subset \pi$ there is $n$ so
that the set of all products of length $n$ of elements of 
$S$ determine fewer than $2^n$ elements of the group. The
formulation we actually use is a 
variation on this.

\begin{lemma} \sl Suppose $\pi$ has subexponential 
growth and $S$ is a 
finite subset of $\pi$. Then the $n$ specified in the 
definition also has the property:
\begin{enumerate}
\item
suppose\/ $T$ is a rooted tree so that all leaves have 
distance $n$ from the root, 
\item each vertex of $T$, other than the leaves and 
the root, has valence $\geq 3$, and
\item each edge of\/ $T$ is labeled with an element of $S$.
\end{enumerate}

Then there is a path in\/ $T$ with distinct endpoints, 
so that the product 
of elements along the path gives $1\in \pi$. 
\end{lemma} 

\proof
First a few clarifications. The
distance of a leaf from the root is the number of vertices 
with valence at least 3 (so branching 
actually occurs) between the leaf and root. Orient each 
edge to go toward the leaves. Then the 
``product of elements along a path'' is the product either 
of the label or its inverse, depending on
whether the path
has the same or opposite direction as the edge.

For each leaf we
get an element  of
$\pi$ by taking the product of elements along the path 
from the root to the leaf. By the distance
condition there are
at least
$2^n$ leaves. The choice of $n$ implies that there are 
two leaves with 
the same associated product. Since
reversing the direction of the path inverts the product, 
the path 
from one of these leaves to the root, then back out
to the other leaf, has total product 1. The geodesic
(unique embedded path) between these leaves also 
has product 1.
Note that the maximum length of 
this path is~$2n$.
\endproof

\Section{Proof of the theorem}

The initial data is a properly immersed grope $G\to M^4$ 
of height at least 2, and 
$\rho\co \pi_1M\to \pi$ with $\pi$
of subexponential growth. We also allow self-intersections 
of the base of the grope, provided their loops have
trivial image in $\pi$. 

Choose dyadic labelings for the branches of $G$, and paths 
from branches to the
basepoint of $M$. Let $S$ be the set 
of images in $\pi$ of $\pi_1$ classes of double point loops in 
$G$, as in~2.4. Let $n$ be the integer associated to
$S$ by the subexponentiality of $\pi$. 

Use Lemmas 2.5 and 2.6 to split $G$ so that it has dyadic 
branches, each branch has a 
$2n$--type, and there are no collisions at
distance $2n$ or less. Note that having $2n$--types 
and having no collisions $\leq 2n$ are both preserved under 
further
dyadic splitting, so doing 
one construction and then the other gives a grope with both 
properties.

Note that if a branch has a cap with no intersections, then we 
can totally contract across that cap. This eliminates
the branch without introducing any new
intersections. Repeating this reduces to the case where all 
caps have an
intersection point.

\proc{Branches with trivial product}
We now claim that if there are any branches 
at all then there is
a path going from branch to branch through at  most
$2n$ intersection points, 
so that the image in $\pi$ of the corresponding loop is trivial. 

Pick some branch, and define a tree 
as follows: At the first stage attach an edge to the root 
corresponding to
each cap of the chosen 
branch. Label the edge by the image in
$\pi$ of the
$\pi_1$ class of the intersection points on the 
cap. Associate to the vertex at the outer end of the edge the
$(2n-1)$--type of the branches 
intersecting the cap. Recall that according to the definition 
of $2n$--type, all of these
branches 
have the same $(2n-1)$--type.

We now do the induction step. Suppose we have gone out 
distance
$k$, with edges associated to intersecting caps and labeled 
by corresponding elements 
of $\pi$, and with vertices
at this distance labeled by
$(2n-k)$--types. Since the edge coming in 
to such a vertex comes from intersecting caps, it 
corresponds to one of the
dyadic cap labels in 
the vertex type. Attach outgoing edges corresponding 
to the other cap labels, label them with
images in $\pi$ of the $\pi_1M$ elements provided by 
the 1--type, and associate to the new vertices 
the $(2n-k-1)$--type
provided by the $(2n-k)$--type and the cap label. The 
number of edges at a 
vertex is the number of caps in
the branch. If the grope has height $\geq 2$ then there 
are at least 
4 caps in each branch. Since the Lemma of
section 3 only requires 3 edges we see that the proof 
works for gropes with height ``3/2''. By this we mean each
branch consists of a cap on one side 
and a grope of height 1 on the other. 

This inductive construction continues until we reach
$k=2n$, 
but in fact we only use the tree out to distance~$n$.

For each branch we now have a tree of 
radius $n$, with root corresponding to this branch, and 
with edges labeled by
elements of
$\pi$.
According to the Lemma of Section 3 there is a path in 
this tree with trivial product. The 
geodesic with the
same endpoints has length
$\leq 2n$ and still has trivial product. 

This verifies the claim of 4.1 in general. There is a 
minor sublety in that this path was constructed 
from abstract
patterns in the type, so we must check that it is 
realized by a sequence of actual
intersections. (Recall that a
type requires that if there are any intersections 
then they satisfy 
certain conditions, but does not  require
that there be any). The types were originally constructed 
to correspond to actual intersections, so the only way this
can fail is that previous steps have
contracted away all the branches that used to lie at 
some point along the path.
However in this 
case  there must be branches with free caps 
somewhere along the path between this point and 
the
root.  This contradicts the standing assumption 
that  all  branches with free caps have already 
been contracted. The
conclusion is that if there are any branches at all 
remaining then paths in 
these trees are realized by actual
intersections.

\proc{Eliminating branches}
The next step is to eliminate the branch at the 
beginning of a path of
the type found in 4.1, without
changing any of the global data ($2n$--types, etc.) 
The path gives us a sequence of
branches$B_i$ 
and dyadic labels
$(x_i,y_i)$ so the
$x_i$ cap in
$B_i$ intersects the $y_i$ cap in $B_{i+1}$. Further,
if $\alpha_i$ is the image in $\pi$ of the loop through this
intersection then the product $\Pi_i\alpha_i$ 
is trivial. This is true for one path of intersections 
starting at
$B_1$. However since $B_1$ has a
$2n$--type it will also be true for {\it all\/} paths 
starting at $B_1$ and following
the same pattern of
dyadic labels. 

Totally contract $B_1$ through the cap $x_1$.  This eliminates $B_1$,
but introduces intersections between the base and the $y_1$ caps of
adjacent branches, and these have $\pi$ images $\alpha_1$. Push these
intersections down and back up through the $x_2$ caps (see 2.2 for the
definition of this operation). This eliminates these intersections,
but introduces intersections between the base and the $y_2$ caps of
branches of distance 2 from $B_1$. These new intersections all have
the same associated $\pi$ elements because $B_1$ had a $2n$--type, and
by construction this element is $\alpha_1\alpha_2$. Note that the
branch $B_2$ (and all other branches that used to intersect the cap
$x_1$) now has a free cap $y_1$, which does not intersect anything
else.  Totally contract these branches along $y_2$, thus reducing the
genus of the base without introducing any new intersections.

Continue in this fashion,pushing $y_i$ intersections down and back up
across $x_{i+1}$, thereby introducing intersections with $y_{i+1}$
caps at distance $i+1$.  As soon as the cap $y_i$ frees up, totally
contract the branch $B_{i+1}$ along it.  At the end of the path all
the base--cap intersections will have the same image in $\pi$, which
by choice of the path is~1. Push these intersections down off branches
to the base. They now are base--base intersections with loop image~1.

We must check that this construction does not introduce any new
base--base intersections until the base--cap intersections are pushed
down at the end. The first step, total contraction of $B_1$ across a
cap, gives an embedded disk because caps have no
self-intersections. Inductively suppose no new base--base
intersections have occurred in step $k$. In each individual branch
pushing down and back up does not introduce intersections. This is
because the transverse sphere used for the operation is embedded.
This in turn follows from the requirement that the path be a geodesic
in the tree, so it enters and exits each branch through different
caps, and there are no self-intersections in the cap it exits
through. Thus the only way intersections can arise is if exit caps on
different branches intersect.  This, however, gives a collision in the
sense of~2.6.  All collisions of distance up to the length of the path
were eliminated by splitting, so they stay disjoint.

Eliminating a branch in this way eliminates 
all intersections with its caps. After doing this we 
check to see if
this leaves branches with caps 
without intersections. If so contract these repeatedly 
until none are left. If there
are any branches 
left after this then repeat the main step (4.1, 4.2). 
Eventually all branches are eliminated, leaving 
a surface satisfying the conclusions of the theorem.
\endproof

\newpage

\centerline{\Large\bf Appendix: Clarification of linear grope height 
raising}

\medskip

\centerline{\sc M Freedman}
\smallskip
\centerline{\sc P Teichner}

\bigskip

Slava Krushkal and Frank Quinn recently brought to our attention
misstatements in the 
proof of our linear grope height raising
procedure which we published in 1995 \cite{FT}. This 
appendix replaces
pages 518--522 of that paper with a proof along the same lines but
with 
correct details. The main difference is that we are more
careful in which order we add surface 
stages. 
This resolves in particular the problem of how to deal with
intersections that involve a dual pair of circles on a surface
stage: Even though the ``key point'' in the middle of page 521 is
not true as stated (the Borromean rings are not slice after all),
the intersections that arise can be dealt with by picking an order
and correspondingly decreasing the scale of the relevant lollipops.  

 We also reformulate
the final word length count in terms of coarse geometry, mainly for
clarity but also for possible future use.

Since the ``warm up'' and ``warm down'' parts of the 
proof of Theorem
2.1 in \cite{FT} are correct,  it suffices to explain the core
construction and show 
that the word length grows linearly. More
precisely, we prove the asserted estimate for the 
word length
$$
\ell(g^\bullet_{k+r})\leq 2r+1 \leqno{(*)}
$$
in terms of the double point loops of 
$G_k$. In the last paragraph on
page 522 this assertion is correctly used to finish the proof of
Theorem 2.1. We now begin the revision on the top of page 518:

\noindent
As we start the 
core construction we have a
Capped Grope $G^c:=G^c_k$ of height
$k\geq 3$. The inductive 
set up is a Grope
$G_{h-1}$ of height $h-1\geq k$ and an embedding
$(G_{h-1},\gamma)\hookrightarrow (G^c,
\gamma)$.  One works with the spines,  proceeding 
from
$g_{h-1}$ to
$g_{h}$ by adding a finite number of connected surfaces $ \Sigma(t)$
to
$g_{h-1}$.
To underline the importance of the order in which the surfaces $
\Sigma(t)$ are 
attached, we write
$$
g_{h-1}=:g(0) \subset g(1) \subset g(2) \subset \dots \subset
g(n)=g_h
$$
where $g(t):=g(t-1)\cup \Sigma(t)$. Even though technically the
$g(t)$ are not gropes 
(since they have heights in between
$h-1$ and
$h$), we will still consider them as such. 
In particular, each $g(t)$
will be thickened to a ``Grope'' $G(t)$. The surfaces
$
\Sigma(t)$ 
are obtained in two steps:

\begin{itemize}

\item Step 1 finds
surfaces $ \Sigma'(t)$ which 
have (illegal)
self-intersections and intersections with grope stages at
various heights, but only 
above
$$
Y := \textrm{base stage} \cup \textrm{second stage
surfaces } 
\Sigma_1\cup \{\Sigma_2\} \textrm{ of } G.
$$

\noindent
The subspace $Y$ is
protected in the construction so that the dual spheres
$\{S\}$ 
will remain geometrically dual to $\{\Sigma_2\}$, the second
stages of $G$, and disjoint from 
everything else.

\item Step 2 only changes the
surface $ \Sigma'(t)$ to $ \Sigma(t)$,
removing double points
  with itself and with earlier stages (and in the process
increases the 
genus of the surface).

\end{itemize}

Every application of Step 1 involves choosing some obvious 
surface
(often a disk) so, formally, the presence of these obvious surfaces
is an inductive 
hypothesis which must be propagated in passing from
$g_{h-1}$ to
$g_{h}$.  The surfaces 
$ \Sigma'(t)$ for Step~1 are of three types:

\begin{enumerate}

\item  ``parallel'' copies of the initial caps $g^c\smallsetminus g$,

\item  meridional disks to some
surface stages of $g(t-1)$, and

\item ``parallel'' copies of stages of the original Grope
$G$.

\end{enumerate}

Every application of Step 2
 is  accomplished by a finite number of
moves called a {\it lollipop move} or a {\it double lollipop move}.
The Step
2 algorithm removes all self-intersections and intersections
of $ \Sigma'(t)$ 
(in a particular order) to produce the surface
$\Sigma(t)$.  The caps
$g^c_h\smallsetminus g_h$,
necessary to define $\ell(g_h)$, are
constructed last and in two steps.  The preliminary caps 
cross all
grope stages above
$Y$ (stages $\ge 3$); these are refined to caps disjoint from the
grope using the dual spheres $\{S\}$.

We next explain the central move in our grope height 
raising
procedure. Every surface stage $\Sigma$  in the Grope $G(t-1)$ has a
symplectic 
basis of circles
$\alpha_1, \beta_1, .. \alpha_g, \beta_g$ where
$g $ is the genus of $\Sigma$, 
along which higher surface stages or
caps have been attached.  We consider tori $T_{\alpha_i} , 
i =1,
\dots, g$ which are $\epsilon$ normal circle bundles to $\Sigma$ in
$G(t-1)$ restricted to 
$\alpha_i$ where $\epsilon$ is a small
positive number depending on $\Sigma$.  Notice that all 
these tori
are disjoint.  Suppose $x$ is a double point  with local sheets $S
\subset \Sigma'(t)$ 
and
$S_\beta \subset \Sigma_\beta$, and
that the surface stage or cap $\Sigma_\beta$ is 
attached to
$\Sigma$ along $\beta$.
Symmetrically, if the surface $ \Sigma'(t)$ intersects 
$
\Sigma_\alpha$ then interchange $ \alpha$ and $ \beta$ in the next
paragraphs.

The {\it lollipop move} replaces a
disk neighborhood $S$ of $x$ with a slightly displaced copy of
$T_\alpha$, made by taking normal
$\epsilon$--bundles over a parallel displacement
(depending 
on $x$) of
$\alpha$ in $\Sigma$, boundary connected  summed to $S$
along a tube which is the 
normal $\epsilon/10$--bundle of
$\Sigma_\beta$ in $G(t-1)$  restricted to an arc $\lambda
\subset
\Sigma_\beta$ from
$(T_{\alpha(\textrm{displaced})})
\cap \Sigma_\beta$ to $x$.  Denote the 
lollipop by $L_{\alpha}$.
  It is the punctured torus made by attaching the tube (or {\it
stem}) to
$T_{\alpha(\textrm{displaced})}$, see Figure 2.1 in \cite{FT}.

We are  now ready to describe the core 
construction in detail. Let
$h-1=k$. The very first application of Step 1 simply attaches one cap
of $g^c$ to $g$.  When regarded as a grope
  stage the self-intersections in the cap are
impermissible and thus the cap only gives $ \Sigma'(1)$.

We specify that the
initial application 
of Step 2 removes (in some order) all
intersections of
$\Sigma'(1)$ using lollipop moves. 
This gives $ \Sigma(1)$ and hence
$g(1)$. To obtain $ \Sigma(2)$ one just repeats Step 1 
and Step 2
by starting with the next cap.
Note that now the self-intersections of the second 
cap as well as
the intersections with the first cap have to be removed (in some
order) by 
lollipop moves.
  In the same manner, one constructs all
surfaces $\Sigma(t)$ and hence the 
grope
$g_{k+1}$. Here the scale $\epsilon$ of
the lollipops is getting rapidly smaller so that they 
do not
intersect the previously constructed surface stages. This is where
the order of things is 
relevant.

In subsequent applications of Step 1 we must specify which surfaces
 we choose 
and what the intersections are.  Each $L_{\alpha}$
contains a meridional circle to which we attach
the meridional disk (type~(2) above) and a longitude
$\ell_{\alpha}$ (picked out by the standard 
framing used to thicken
$g$ to $G$) to which we attach a 
surface of type~(1)--(3) above. 
Type~(2) arises if $\ell_{\alpha}$ is the meridional circle of a
previously constructed lollipop. Types~(1) or (3) arise if
$\ell_{\alpha}$ is parallel to a circle on the original grope $g$.
In this case, the new surface
or cap is only crudely 
parallel in the sense that we need to glue an
annulus
$A$ to get from the longitude $\ell_\alpha$ 
to $ \partial
\Sigma_{\alpha(\textrm{displaced})}$, the attaching circle of a
slightly displaced copy 
of one of the surfaces or caps of $G^c$. The
surface
$
\Sigma'(t)$ is then defined to be 
$A\cup
\Sigma_{\alpha(\textrm{displaced})}$.
The framing assumption of $G$ implies that for type~(3) 
the
surface stage
$\Sigma_{\alpha(\textrm{displaced})}$ will be disjoint from
everything constructed
previously, ie from $g(t-1)$. However, for
both types~(1) and (3), the annulus
$A$ may intersect
many $ \Sigma(s), s<t$, so that $
\Sigma'(t)$ has many intersections with $g(t-1)$.
For type~(2), 
$ \Sigma'(t)$ is a meridional disk
and it will intersect $g(t-1)$ in a single point.

The reader may
expect that the next application of Step 2 will use lollipop moves on
$ \Sigma'(t)$ to remove these
intersection points.
This is  part of the picture, but there is a difficulty.  The
lollipop moves, if 
repeated, produce a branch heading inexorably
down $G$: namely resolving (meridian disk)
$\cap\Sigma_i$ with a lollipop capped by a (meridian disk) meeting a
$\Sigma_{i-1}$ lead 
toward the base of $G$ which is $\Sigma_1$.
There is no way of using a lollipop to remove a 
point of (meridian
disk)
$\cap \Sigma_1$.  The solution is to use the {\it double lollipop
move} 
to resolve any intersection of a current top stage meridional
disk with a third stage surface 
$\Sigma_3$.  This move turns the
branch of the growing grope back ``upward" to avoid the 
bottom part
$Y$.

The {\it double lollipop move} removes an intersection
$x$ between a surface
$\Sigma'(t)$  and a third story surface $\Sigma_3$.  This move
replaces a small disk neighborhood 
$S \subset \Sigma'(t)$ of $x$
with $L_\alpha/\Sigma_\alpha$.  The notation assumes
$\Sigma_3$
attaches to $\beta$ (otherwise reverse the labels
$\alpha$ and $\beta$), $L_\alpha$ is the lollipop 
made from
$T_\alpha$ as describe above,
$\Sigma_\alpha$ is the third story surface attached to
$\alpha$ and finally $L_\alpha/\Sigma_\alpha$ denotes the embedded
surface that results by 
surgering
$L_\alpha$ along a parallel copy $\Sigma_{\alpha(\textrm{displaced})}$
of
$\Sigma_\alpha$, 
ie,
$L_\alpha/\Sigma_\alpha = (L_\alpha\smallsetminus
\textrm{nbh. of } 
\alpha(\textrm{displaced})) \cup$ 
two copies of
$\Sigma_{\alpha(\textrm{displaced})}$.  Because we have assumed
$G^c$ is an 
untwisted thickening the two copies of
$\Sigma_{\alpha(\textrm{displaced})}$ are disjoint from  each 
other
and from the original $\Sigma_\alpha$.

Now suppose that we have constructed the grope
$g_{h-1}$. Then the
top layer of surfaces has a natural symplectic basis coming from the
original 
grope $g$ and the (meridian, longitude) pair on each
lollipop. These bound obvious surfaces
$ \Sigma'(t)$ of types~(1)--(3) as explained above.
Applying Step 2 to these surfaces in some 
chosen order, we remove
intersection points by a lollipop move except in the case of
intersection 
with a third stage surface $\Sigma_3$ in which case a
double lollipop is used.
This gives the 
embedded surfaces $ \Sigma(t)$ and hence an embedded
grope
$(g_{h},\gamma)\hookrightarrow 
(G^c, \gamma)$.

We next check the normal
framing. If we assume that each cap has
algebraically 
zero many self-intersections then
all surfaces $ \Sigma'(t)$ are $0$--framed.
A lollipop move on a
$\pm$--self-intersection
changes the relative Euler class by $\pm 2$ (this is best
checked in the
closed case, $S^2\times S^2$, where adding
the framed dual $0 \times S^2$ to
$S^2\times 0$ 
gives the diagonal). All other lollipop moves leave
the $0$--framing unchanged.
Thus the passage 
to $ \Sigma(t)$ leaves the relative Euler class
trivial so the neighborhood of $g(t)$ agrees with 
the standard
thickening $G(t)$.

To obtain caps $\{ \delta\}$ for $g_h$, we examine the symplectic
basis for the top stage of $g_h$. Some of the curves bound meridian
disks to earlier stages of the
construction. Some bound
``parallel'' copies of sub capped gropes of $G^c$. Contracting, the
latter 
also yield disks.
We set $h=k+r$ and
$$g^\bullet_{k+r} := g_{k+r} \cup \{\delta\}$$ 
The superscript $\bullet$ warns the reader that $g^\bullet_{k+r}$ does
not satisfy the definition of a capped grope owing to the cap--grope
intersections.  These will be removed in the last step, see the last
paragraph of page 522 in \cite{FT}.

Let us next bound the word length
$\ell (g^\bullet_{k+r})$ in terms of the original
generators (=
double point loops) of the free group
$F:=\pi_1G^c$. Recall that we need to prove
$$ \ell(g^\bullet_{k+r})\leq 2r+1. \leqno{(*)}
$$

For this purpose, we put a {\it pseudo metric} on the
universal
covering $X$ of $G^c$.
This is a distance function which still satisfies the
triangle 
inequality but distinct points may have distance zero. Note
that pseudo-metrics can be pulled back 
by arbitrary maps which we
will use in the construction as follows.
First project
$X$ onto the Cayley
graph of $F$ such that lifts of the Grope body
$G$ map bijectively onto the vertices and lifts of the
plumbed
squares in the Caps map bijectively onto the centers of the edges.
Then take a coarse or
pseudo version of the usual path metric on the
Cayley graph (in which all edges have length~1) by
saying that edge
centers have distance~1/2 from all the vertices the edge meets and
that all path
components of the Cayley graph minus the edge centers
have diameter zero.
Finally, use the above 
map to pull this pseudo metric back to~$X$.

For any map $f\co Y\to G^c$ which is trivial on $\pi_1$, we
may then
measure the {\it diameter} of a lift $\tilde f(Y)$ in $X$. For
example, if $Y$ is a model 
capped grope (ie with unplumbed caps)
such that $f(Y)=g^\bullet_{k+r}$ then the diameter of 
$\tilde f(Y)$
is just the word length $\ell(g^\bullet_{k+r})$.

If $Y$ happens to be a disk, surface 
or (capped) grope such that $
\partial Y$ maps to $G$, it is very useful to consider the {\it
radius} 
of $\tilde f(Y)$ around the ``point'' $\tilde f( \partial
Y)$. This uses the fact that each lift of $G$ 
projects onto a vertex
in the Cayley graph of $F$ and thus has radius zero itself. For
example, if $Y$ 
is a disk mapping onto a cap of $G^c$ which has one
self-intersection, then the radius of $\tilde f(Y)$
is~$1/2$ whereas
the diameter is~$1$.

Let $X_r$ be a lift of $g_{k+r}$ to $X$ and let $X_r^c:= \tilde
f(Y)$ where $f(Y)=g^\bullet_{k+r}$ as above. Then the triangle
inequality shows that 
$radius(X_r^c)\leq radius(X_r)+1/2$ and hence
$$
\ell(g^\bullet_{k+r}) = diam(X_r^c)\leq 2\cdot radius(X_r^c)\leq
2\cdot radius(X_r)+1.
$$
It thus suffices to check that $radius(X_r)\leq r$.
This in turn follows by the triangle inequality (applied to
the usual tree structure of the grope 
$g_{k+r}$) from knowing that
the radii of all
$S(t)$ are $\leq 1$. Here
$S(t)$ are lifts to
$X$ of the
surfaces $\Sigma(t)$ used in the construction of
$g_{k+r}$ and the radii are again measured w.r.t. 
$ \partial S(t)$.

We prove that radius~$S(t)\leq 1$ by induction on $t$: Recall that
the first surface 
$ \Sigma(1)$ was obtained by applying lollipop
moves to the first cap of $G^c$. Before the lollipop
moves, we can
lift the (unplumbed) cap to $X$ and as explained above it has
radius~$1/2$ (if
the 
cap is embedded then the radius is zero but we won't
consider this easy case). The lollipops then
increase this radius
to at most~$1$, independently of how many are used. This follows from
the 
triangle inequality applied to the decomposition of each lollipop
into its stem and body (or toral piece).
The body has diameter zero
since it lies in $G$ whose lift projects to a vertex. The stem has
by
definition diameter~$1/2$ since it leads from a plumbed square to
the base of the cap.

Now assume by induction that radius~$ S(s)\leq 1$ for all $s<t$.
Let $S'(t)$ be a lift to $X$ of 
$ \Sigma'(t)$.
If $ \Sigma'(t)$ is of type~(2) or~(3) then the radius of $S'(t)$ is
zero since it lies 
in a lift of $G$. For every intersection point of
$\Sigma'(t)$ with $g(t-1)$ we add a lollipop or a 
double lollipop
to obtain $
\Sigma(t)$. Only the stems of these (double) lollipops will
contribute 
to the radius of
$S(t)$ since the bodies lie in $G$.
The induction hypothesis implies that all these 
stems have diameter
$\leq 1$ and thus we are done in this case.

Finally, consider the case where
$\Sigma'(t)$ has type~(1), ie is a ``parallel'' cap. Then its
radius is $1/2$ as explained above. 
For every self-intersection of
$\Sigma'(t)$ and every intersection point of
$\Sigma'(t)$ with 
$g(t-1)$ we add a lollipop
to obtain $\Sigma(t)$ (note that double lollipops don't
occur for caps). 
Again, only the stems of these lollipops
will contribute to the radius of $S(t)$.
There are two types 
of lollipops:
One type removes self-intersections and intersections with surface
stages of $g(t-1)$ 
that come from the caps of $g^c$.
  As for $\Sigma(1)$ the corresponding lollipop stems have
diameter~$1/2$ and thus can only increase the radius to~$1$. The
other type of lollipops remove
intersections of the annulus
$A=($collar of $ \partial
\Sigma'(t))$. This means that, as far as our 
pseudo metric can
measure, the stems of the lollipops start essentially on $ \partial
\Sigma'(t)$ 
which is the base point with respect to which we measure
the radius. By the induction hypothesis 
these stems can only bring
the radius up to~$1$.
\endproof

{\bf Note added in proof}\qua Slava Krushkal has pointed out that in
the above proof, the ``warm-up'' and ``warm-down'' steps can be
replaced by the following easier and shorter argument:

Do the core
construction on the originally given Capped Grope of 
height $k\geq 2$, preserving only the bottom
surface $\Sigma_1$ 
instead of the first two stages $Y$ as done above. (No dual spheres 
need to be
constructed.)
After the core construction, we have a Capped Grope of height~$k+r$ 
and word length
$\leq 2r+1$, with
many cap--body intersections but caps are disjoint from the
bottom surface
$\Sigma_1$. Now do symmetric contraction
of the bottom surface. This requires taking parallel
copies
of whatever is attached to it, and reduces the height
of the entire Capped Grope by~1. Then push all
cap--body intersections 
down and off the contraction. This at most doubles
the estimate on the double
point loop length and thus leads to a 
clean Capped Grope of height~$k+(r-1)$ and word length
$$ \leq 2(2r+1)=4(r-1)+6.
$$
Thus linear grope height raising is established.

\end{document}